\documentclass[invmat,final]{svjour}
\usepackage{amsfonts}
\newcommand{\QED}{{\unskip\nobreak\hfil\penalty50%
\hskip1em\hbox{}\nobreak\hfil $\qed$%
\parfillskip=0pt \finalhyphendemerits=0 \par\medskip\noindent}}
\newcommand{\n}{\par\noindent}

\newcommand{\sn}{\par\smallskip\noindent}
\newcommand{\pars}{\par\smallskip}
\newcommand{\parb}{\par\bigskip}
\newcommand{\euK}{\mathfrak{K}}
\newcommand{\euL}{\mathfrak{L}}
\begin{document}
\title{A correction to Epp's paper ``Elimination of
wild ramification''}
\author{Franz-Viktor Kuhlmann
}                     
%
%
\institute{
 Mathematical Sciences Group,
 University of Saskatchewan,
 106 Wiggins Road,
 Saskatoon, Saskatchewan, Canada S7N 5E6
}
\date{}
%
\maketitle
\begin{abstract}
We fill a gap in the proof of one of the central theorems in Epp's
paper, concerning $p$-cyclic extensions of complete
discrete valuation rings.
\end{abstract}
%
%
In his famous paper \cite{1}, Epp considers the following situation: $S$
and $R$ are two discrete valuation rings such that (1) $S$ dominates
$R$, and (2) if the characteristic $p$ of the residue field of $S$ is
not zero, then its largest perfect subfield is separable and algebraic
over the residue field of $R$. He proves that then, there exists a
discrete valuation ring $T$ which is a finite extension of $R$ such that
the localizations of the normalized join of $S$ and $T$ are weakly
unramified over $T$. Towards this result, he proves the following
theorem, assuming that all discrete valuation rings are complete:
\sn
{\bf Theorem~(1.3).} \ {\it Let $S$ be a $p$-cyclic extension of $S_0$
where $S_0$ is a weakly unramified extension of $R$ such that
${\euL}^{p^{\infty}}={\euK}$, where $\euL$ and $\euK$ are the residue
fields of $S_0$ and $R$ respectively. There exists a finite extension
$T$ of $R$ such that $TS$ is weakly unramified over $T$.}

\pars
There is a mistake in the proof of the {\it Equal characteristic
$p\ne 0$ case}. We will sketch those parts of the proof that are
necessary to understand and correct the mistake.

\pars
Using well-known structure theorems that are discussed in section 0.1 of
his paper, Epp writes $R={\euK}[[\pi]]$ and $S_0={\euL}[[\pi]]$, where
$\pi$ is a local parameter of $R$. By Artin-Schreier theory and the fact
that the Artin-Schreier polynomial $X^p-X$ is additive and surjective on
the maximal ideal of the power series ring $S_0\,$, Epp finds that the
$p$-cyclic extension $S$ of $S_0$ is defined by an equation of the form
\begin{equation}                            \label{de}
z^p-z\>=\>a_{-N}\pi^{-N}+\ldots+a_{-1}\pi^{-1}+a_0
\end{equation}
with $a_i\in {\euK}$. In the case of $N=0$ there is nothing to prove,
so we assume that $N\ne 0$. Epp defines the following subsets of
$\{1,\ldots,N\}$:
\[I\>=\>\{m\mid a_{-m}\in {\euK}, a_{-m}\ne 0\}\;\mbox{ \ \ and \ \ }
\; J\>=\>\{m\mid a_{-m}\notin {\euK}, a_{-m}\ne 0\}\;.\]
After dealing with the case of $J=\emptyset$, Epp assumes that
$J\ne\emptyset$. The idea is now to find some $d\in {\euK}[[t]]$, where
$t^{p^k}=\pi$ for some $k$, such that after replacing $z$ by $z-d$
and adding $d-d^p$ to the right hand side of equation (\ref{de}), the
defining equation (of ${\euK}[[t]]S$ over ${\euK}[[t]]$) will be of
the form
\[z^p-z=ct^{-n}+\ldots\]
where $n$ is divisible by $p$, but $c\notin {\euL}^p$. Epp shows that
then, ${\euK}[[t]]S$ is weakly unramified over ${\euK}[[t]]$, so we
can take $T={\euK}[[t]]$.

Note that a transformation of the above type replaces a $p$-th power
$d^p$ on the right hand side of (\ref{de}) by its $p$-th root $d$. By a
repeated application of such replacements, Epp seeks to get rid of all
coefficients that lie in ${\euL}^p\setminus {\euK}$. To this end, he
chooses a positive integer $\nu$ such that
\[\min_{m\in J} m p^\nu \> > \> \max_{m\in I} m\;.\]
Since ${\euL}^{p^{\infty}}={\euK}$,
\[\nu_m\>:=\>\max \{i\mid a_{-m}\in {\euL}^{p^i}\}\><\>\infty\]
for every $m\in J$. Let $\mu:=\max\{\nu_m\mid m\in J\}$, and let $t$
be such that $t^{p^{\nu+\mu}}=\pi$. Equation~(\ref{de}) can now be
written
\[z^p-z\>=\>\sum_{m\in J}a_{-m} t^{-mp^{\nu+\mu}}\>+\>
\sum_{s\in I}a_{-s} t^{-sp^{\nu+\mu}}\>+\> a_0\;.\]
Using the above described transformations, Epp arrives at a defining
equation
\begin{equation}                            \label{de1}
z^p-z\>=\>\sum_{m\in J}c_{-m} t^{-mp^{\nu+\mu-\nu_m}}\>+\>
\sum_{s\in I}c_{-s} t^{-s}\>+\> a_0\;,
\end{equation}
where:\n
$\bullet$ \ for every $s\in I$, $c_{-s}\in {\euK}$ is the
$p^{\nu+\mu}$-th root of $a_{-s}\in {\euK}$ (note that ${\euK}$ is
perfect!),\n
$\bullet$ \ for every $m\in J$, $c_{-m}\in {\euL}\setminus {\euL}^p$
is the $p^{\nu_m}$-th root of $a_{-m}$ (recall that $a_{-m}\in
{\euL}^{p^{\nu_m}}\setminus {\euL}^{p^{\nu_m+1}}$),\n
$\bullet$ \ for every $m\in J$, $p$ divides $mp^{\nu+\mu-\nu_m}$ since
$\nu\geq 1$ and $\mu\geq \nu_m\,$,\n
$\bullet$ \ for every $m\in J$ and $s\in I$, $-s>-mp^{\nu+\mu-\nu_m}$
since $s< m p^\nu$ by the choice of $\nu$.

\pars
Now Epp claims that the term in $t$ with the most negative exponent has
a coefficient which is not in ${\euL}^p$. This is not necessarily
true. It would hold if the exponents $-mp^{\nu+\mu-\nu_m}$ were
distinct, for distinct $m$. But this could be false since we know
nothing about the $\nu_m\,$.

\sn
{\bf Example.} \ Suppose that
\[z^p-z\>=\>a_{-m_0}\pi^{-m_0}+a_{-m_1}\pi^{-m_1}+\ldots\]
where $m_0=m_1p$, $a_{-m_0}=c_0^p$ and $a_{-m_1}=c_1-c_0$ with $c_0\in
{\euL} \setminus {\euL}^p$ and $c_1\in {\euL}^p$. Then $a_{-m_0}
\notin {\euK}$, $a_{-m_1} \notin {\euL}^p$, $\nu_{m_1}=0$,
$\nu_{m_0}=1$, and $a_{-m_1}+ a_{-m_0}^{1/p}= a_{-m_1}+c_0=c_1\in
{\euL}^p$. Using the notation of (\ref{de1}), we find that
$-m_0p^{\nu+\mu-\nu_{m_0}}= -m_1p^{\nu+\mu-\nu_{m_1}}$ and
$c_{-m_0}+c_{-m_1}=c_1\in {\euL}^p$.                              \QED

\pars
So we see that the coefficient of the term in $t$ with the most negative
exponent can well lie in ${\euL}^p$. Choosing $c_1$ to lie in ${\euK}$
in our example, we see that this coefficient may even lie in ${\euK}$,
so that the corresponding exponent ``switches'' from the set $J$ to the
set $I$. However, whenever such a recombination happens and we start
over with the new equation (\ref{de1}), the new set $J$ will be smaller
than the original set $J$. So the gap in Epp's proof can be closed by
repeating his transformations until his assertion is satisfied or the
set $J$ is empty.

The latter may well happen: consider the equation
\[z^p-z\>=\>a_{-m_0}\pi^{-m_0}+a_{-m_1}\pi^{-m_1}+a_0\]
with the conditions of our example, and assume in addition that $c_1=0$.
Then the transformation leads to the equation
\[z^p-z\>=\>a_0\;.\]
This shows that even if $J\ne \emptyset$ in the original equation
(\ref{de}), the residue field extension of $TS$ over $T$ may end up to
be an Artin-Schreier extension, in contrast to the purely inseparable
extension which Epp obtains for this case.

\parb
A far-reaching generalization of Epp's results will be proved in
\cite{2}.


\begin{thebibliography}{}
\bibitem{1}
Epp, Helmut P.: Eliminating Wild Ramification. Invent. Math.
\textbf{19} (1973) 235--249
\bibitem{2}
Kuhlmann, F.--V.: The generalized stability theorem
and henselian rationality of valued function fields. In preparation
\end{thebibliography}
\end{document}